\documentstyle[11pt,righttag,verbatim]{amsart}
\newtheorem{thm}{Theorem}[section]

\newtheorem{prop}{Proposition}[section]

\theoremstyle{definition}

\theoremstyle{remark}

\def \mrn {{\Bbb R}^n}

\def \eps {\epsilon}

\def \p {\partial}

\def \rao#1 {\frac{\p}{\p #1} #1}

\def \la {\langle}
\def \ra {\rangle}
\def \s {S}

\numberwithin{equation}{section}

  \title{Magnetic Inverse Problem}
   \author {Mark S. Joshi}
   \address{Department of Pure Mathematics and Mathematical
Statistics, University of Cambridge, 16 Mill Lane, Cambridge CB2 1SB,
England, U.K.} 
   \email{joshi@@dpmms.cam.ac.uk}
\author{Antonio Sa Barreto}
   
\begin{document}
\section{Computing the Symbol}

In this section, we compute the symbol of the scattering matrix - in
particular given two magnetic fields which agree to some order, we
compute the principal symbol of the difference of the associated
scattering matrices in terms of the lead term of the difference of the
magnetic fields.  We use the techniques of \cite{smagai},
\cite{smaginv} and \cite{explicit} to construct the Poisson operator
and compute the symbols.  We proceed explicitly where
possible but occasionally fall back on results from \cite{smagai} for
brevity. 

Let $P$ be the operator,
\begin{equation}
P =  \sum \limits_{j=1}^{n} \left( i \frac{\p}{\p x_j} + A_j \right)^2
+ V  = \Delta + i \sum \limits_{j=1}^{n} A_j \frac{\p}{\p x_j} + q
\end{equation}
where $(A_1,\dots,A_n)$ and $V$ are real-valued classical symbols of order
$-2.$ The Poisson operator is then the map that maps $f \in
C^{\infty}(S^{n-1})$ to the smooth function $u$ such that $(P -
\lambda^2)u=0,$ and 
$$u = e^{i\lambda |x|}|x|^{-\frac{n-1}{2}} f(x/|x|) + 
e^{-i\lambda |x|}|x|^{-\frac{n-1}{2}} g(x/|x|) +
O(|x|^{-\frac{n+1}{2}}).
$$ The kernel of the Poisson operator will be a smooth function on
$S^{n-1} \times \mrn$ with singular asymptotics. The scattering
matrix is the map on $C^{\infty}(S^{n-1}),$ 
\begin{equation}
S(\lambda) : f \mapsto g.
\end{equation}

Our first result in this section is
\begin{prop}
With $(A_1, \dots, A_n)$ and $q$ as above, we have that $a^{*}S(\lambda)$
is a zeroth order classical pseudo-differential operator, where
$a(\omega) = -\omega.$
\end{prop}
We remark that saying that $a^{*}S(\lambda)$ is a zeroth order
classical pseudo-differential operator is equivalent to saying that is
a classical Fourier integral operator of order $0$ associated with
geodesic flow at time $\pi.$ This proposition is clear from any of
\cite{smagai}, \cite{explicit} and \cite{andras} by just observing that
the arguments are not changed by adding a first order short range
self-adjoint perturbation, we therefore only present the parts of the
proof which relate to the proof of our main result:

\begin{thm} \label{symbolcomputation}
Let $(A_1,\dots,A_n)$ and $(A^{'}_{1},\dots,A^{'}_{n})$ be real-valued
classical symbols of order $-2$ with $S(\lambda)$ and $S'(\lambda)$
the associated scattering matrices. 
 If $A_j - A^{'}_{j}$ is of
order $-k$ with  $k \geq 2$ then $S(\lambda) - S'(\lambda)$ is of
order $1-k.$  If $\sum \limits_{j=1}^{n} (A_j - A^{'}_{j}) dx_j = B$
with lead term $B^{(k)},$ and $B^{(k)}$ is aradial then the principal
symbol of $S(\lambda) -
S'(\lambda)$  determines and is determined by
$$\int \limits_{0}^{\pi} \langle B^{(k)}(\gamma(s)), \frac{d\gamma}{ds}(s)
\rangle (\sin s)^{k-1} ds$$
for all geodesics $\gamma,$ where we regard $B^{(k)}= \sum B^{(k)}_{j}
dx_j$ as a one form canonically pairing with the vector
$\frac{d\gamma}{ds}(s).$ 
\end{thm}
We say a one-form is {\it aradial} if its pairing with the radial vector
field, $x\frac{\p}{\p x},$ is zero. 

Following the ideas of \cite{smagai}, \cite{explicit},
\cite{smaginv}, we look to construct the Poisson operator for the
problem as a sum of oscillatory integrals and then use this to read
off the properties of the scattering matrix. 
In particular, we attempt to construct the Poisson operator as 
\begin{equation}
e^{i\lambda x.\omega}(1+b(x,\omega)),
\end{equation}
with $b$ a classical symbol of order $-1$ in $x$ and smooth in $\omega.$
We will see that this ansatz works away from the set $\omega = - x/|x| .$
Applying $P-\lambda^2,$  we obtain
$$e^{i\lambda x.\omega} \left(  -2i\lambda \omega.\frac{\p b}{\p x} + 
\Delta b - \lambda \sum \limits_{j=1}^{n} \omega_j A_j (1+b)  + i
\sum \limits_{j=1}^{n} A_j \frac{\p b}{\p x_j} +q +qb
\right).$$ 
That this can be solved smoothly to infinite order in a neighbourhood of
$x/|x|=\omega$ is just a repetition of the argument in the proof of
Proposition 18 of \cite{smagai}.  

We let $b$ have asymptotic expansion $\sum \limits_{j=1}^{\infty}
b_{-j},$ let $A_{l}$ have expansion $\sum \limits_{j=2}^{\infty} A_{l}^{(-j)}$ 
and $q$ have expansion $\sum  \limits_{j=2}^{\infty} q_{-j}.$ 
In order to continue $b$ smoothly we observe that the lead term is 
$$-2i\lambda \omega.\frac{\p b}{\p x} - \lambda \sum \limits_{j=1}^{n}
\omega_j A^{(-2)}_j +q_{-2}.$$ We want this to be zero. 

Note as these are homogeneous functions this is really an equation on the
sphere. We therefore take coordinates $(r,s,\theta)$ where $r=|x|,$ $s$ is the 
geodesic distance of $x/|x|$ from $\omega$ and $\theta$ is the angular
coordinate about $\omega.$ Note the coordinate system depends on $\omega$ but
we shall suppress $\omega$ in our notation most of the time. 

Now with out loss of generality,
 we can take $\omega$ to be the north pole. Then 
$$ \omega . \p_x = \p_{x_n} .$$ Now $\cos(s) = \frac{x_n}{|x|}, $ $r=|x|.$ The
$\theta$ coordinate will be purely parametric. We have
$$\frac{\p s}{\p x_n} = \frac{-1}{\sqrt{1- \frac{x_{n}^{2}}{|x|^2}}} \left(
\frac{1}{|x|} - \frac{x_{n}^{2}}{|x|^3} \right) = -r^{-1} \sin(s),$$ 
and $$\frac{\p r}{\p x_{n}} = \cos(s).$$ 
So applying $\omega.\p_{x}$ to $b_{-j}(s,\theta; \omega)r^{-j},$ we obtain
$$r^{-j-1} \left[ -\sin(s) \frac{\p b_{-j}}{\p s} - j \cos(s) \right] b_{-j}.$$ 
Let $W_{-2} = -\lambda \sum \limits_{j=1}^{n}
\omega_j A^{(-2)}_j - q_{-2}.$ 
So for $j=1,$ we have taking $r=1$
$$ 2i\lambda (\sin(s) \p_{s} + \cos(s) ) a_{-1} = W_{-2},$$ which is equivalent
to 
$$ 2i\lambda \p_s (\sin(s) a_{-1} )  =   W_{-2}.$$ 
We want $b$ to be smooth at $s=0,$  so
$$\sin(s) b_{-1} =  \frac{1}{2i\lambda} \int \limits_{0}^{s}  W_{-2}  ds',$$
which implies that
\begin{equation} a_{-1}(s,\theta; \omega) = \frac{i}{2\lambda \sin(s)} \int \limits_{0}^{s}
W_{-2}(s',\theta; \omega) ds'. \label{firsterrorterm}
\end{equation}
This will be singular as $s \to \pi- ,$ but let's ignore that for now. Now
suppose we have chosen the first $j$ terms so we have an error $d \in
S^{-j-1}_{cl} $ with lead term $d_{-j-1}(s,\theta; \omega)r^{-j-1}.$ We then want to
solve the transport equation,
$$-2i\lambda \omega . \p_z ( a_j r^{-j}) + d_{-j-1} r^{-j-1} =0,$$
as above we get 
$$  2i \lambda [ \sin(s) \p_s +j \cos(s) ] a_j +d_{-j-1} =0.$$
We solve this to obtain,
\begin{equation}
a_j(s,\theta, \omega) = \frac{i}{2\lambda (\sin(s))^{j}} \int \limits_{0}^{s} 
(\sin(s'))^{j-1} d_{-j-1}(s', \theta; \omega) ds'. \label{genterm}
\end{equation}
So away from $s=\pi,$ we can achieve an error in $S^{-\infty}$ by applying
Borel's lemma. ie away from the antipodal point. Note that we have a focussing
of the geodesics and as well as the fact the solutions blow-up we also have
that they will have different values according to the angle. In
particular provided $d_{-j-1}$ does not grow faster than $(\pi
-s)^{1-j}$ we have that $a_j$ grows as $(\pi -s)^{-j}.$ 

Before introducing a second ansatz to cope with the antipodal point, 
we compare the Poisson operators associated to two different magnetic
potentials. Suppose $(A_1, \dots, A_n)$ and $(A^{'}_1, \dots,
A^{'}_n)$ are both classical symbols of order $-2$ and the difference
is $(B_1, \dots, B_n)$ which is a classical symbol of order $-k,$ with
lead term $(B^{(-k)}_{1}, \dots, B^{(-k)}_{n}).$ 
The
first $k-2$ forcing terms above are then unchanged and the forcing
terms at level $k-1$ will differ by $W_{-k} = -\lambda \sum \limits_{j=1}^{n}
\omega_j B^{(-k)}_{j} .$ (Note that the change in the zeroth order term
will be lower order.)

Thus the lead term of the difference of the Poisson operators will be
\begin{equation}
\frac{i r^{1-k}}{2\lambda (\sin s)^{k-1}} \int \limits_{0}^{s}
W_{-k}(s',\theta; \omega) (\sin s')^{k-2} ds'.
\end{equation}
This is the important result in our construction as we will see that
the lead term of this as $s \to \pi-$ is essentially the principal symbol of
the difference of the scattering matrices. We therefore want an
invariant interpretation of 
$$\int \limits_{0}^{\pi} (\sin s)^{k-2} W_{-k}(s,\theta';\omega) ds.$$
If we take $\omega$ to be the north pole and rotate so that $\theta' =
(1,0,\dots,0),$ the computation lies entirely in the $(x_1,x_n)$
plane and $W_k(s)$ equals $-\lambda B_{n}^{(-k)}(s).$ Now if we assume
$B$ is aradial then an elementary computation shows that,
$$-\langle B^{(-k)}(\gamma(s)), \frac{d\gamma}{ds}(s) \rangle \sin s=
B_{n}^{(-k)}$$ in this case. So by rotational invariance we deduce
that in general the lead term of the difference of the Poisson
operators is 
\begin{equation}
\frac{i r^{1-k}}{2(\sin s)^{k-1}} \int \limits_{0}^{s} \langle
B^{(-k)}(\gamma(s')), \frac{d\gamma}{ds'}(s') \rangle (\sin s')^{k-1}
ds'
\end{equation}
and that the lead singularity as $s \to \pi-$ is
\begin{equation}
\frac{i r^{1-k}}{2(\pi-s)^{k-1}} \int \limits_{0}^{\pi} \langle
B^{(-k)}(\gamma(s')), \frac{d\gamma}{ds'}(s') \rangle (\sin s')^{k-1}
ds'.
\end{equation}

The remainder of the construction of the Poisson operator and the
computation of the symbol is now just
a repetition of the arguments in \cite{smagai} or \cite{explicit}. 
We sketch these for completeness.

Taking $\omega$ to be the north pole, close to the south pole we look
for the Poisson operator in the form, 
So close to the south pole, we look for the Poisson operator in the form,
\begin{equation}
\int \limits_{0}^{\infty} \int \fracwithdelims(){1}{\s|x|}^{\gamma} \s^{\alpha}
\fracwithdelims(){1}{S|x|}^{\gamma} S^{\alpha} e^{i(Sx'.\mu -
\sqrt{1+S^2}|x|)} a\left(\frac{1}{S|x|},S,\mu \right) dS d\mu,
\label{secansatz} 
\end{equation}
with $a(t,\s,\mu)$
 a smooth function compactly supported on $[0,\eps) \times [0,\eps)
\times S^{n-2}$ and $\alpha = \frac{n-3}{2}, \gamma = -\frac{n-1}{2}.$
We assume that $\omega$ has been rotated to the north pole.
Note that for $|x|$ in a compact set the integral is supported on a compact set
and so we have no problems with convergence - in particular the integral yields
a smooth function.

In the lower hemi-sphere, away from the south pole, this ansatz is
equivalent to the original one - this follows from an application of
stationary phase  (see \cite{explicit}). 
However the lead term in $|x|$ at order $-k$ is now
allowed to be singular of order $-k$ as $s \to \pi-$ (which
corresponds to $S=0+$) and there is no
constraint on the values for different angles matching. This allows
the transport equations to be solved right up to the antipodal point
and to all orders. The error is then removed by applying the resolvent
which yields a term of the form
$e^{-i\lambda|x|}|x|^{-\frac{n-1}{2}}h(x)$ with $h$ a classical zeroth
order symbol - this term will not affect the singularities in the
distributional asymptotics of the Poisson operator. 

We recall Proposition 3.4 of  \cite{explicit}, which is a special case
of Proposition 16 of \cite{smagai}. 
\begin{prop} \label{disasymp}
If $u(x,\omega)$is of the form \ref{secansatz} then 
$ e^{i\lambda|x|} \int u(|x| \theta,\omega) f(\theta, \omega) d\theta
d\omega$ is a 
smooth symbolic 
function in $|x|$ of order $-1-\alpha$ and its lead coefficient is
$|x|^{-\alpha-1} \la K, f
\ra$ where $K$ is the pull-back of the Schwartz kernel of a
pseudo-differential operator of order $\alpha-\gamma-(n-2)$ by the map
$\theta \mapsto
-\theta.$ The principal symbol of $K$ determines and is determined by the lead
term of the symbol, $a(t,\s,\mu),$ of $u$ as $\s \to 0+.$
\end{prop}

The fact that the scattering matrix is the pull-back of a
pseudo-differential operator is now immmediate. To deduce Theorem
\ref{symbolcomputation}, we observe that from our computations above
we have that the difference of the second ansatzs for the two Poisson
operators will be of the same form but with $\alpha$ increased by
$k-1$ and the lead term of the symbol as $S \to 0+$ will be a constant
multiple of 
$$S^{n-k-1}  \int \limits_{0}^{\pi} \langle
B^{(-k)}(\gamma(s')), \frac{d\gamma}{ds'}(s') \rangle (\sin s')^{k-1}
ds'.$$ 
So the theorem then follows from Proposition \ref{disasymp}.


\begin{thebibliography}{aa}
\bibitem{eskral} G. Eskin, J. Ralston, Inverse Scattering Problem for
the Schrodinger Equation with Magnetic Potential at a Fixed Energy,
Commun. Math. Phys. 173, 199-224 (1995)
\bibitem{helg} S. Helgasson, {\em Groups and Geometric Analysis,}{
Academic Press 1984.}
\bibitem{coulomb} M.S. Joshi, {\em Recovering Asymptotics of
Coulomb-like Potentials,} to appear in S.I.A.M. Journal of
Mathematical Analysis
\bibitem{explicit} M.S. Joshi, {\em Explicitly Recovering Asymptotics
of Short Range Potentials,} preprint
\bibitem{smaginv} M.S. Joshi, A. S\'{a} Barreto, {\em Recovering Asymptotics of Short
Range Potentials,} Comm. Math. Phys. 193, 197-208 (1998)
\bibitem{metric} M.S Joshi, A. S\'{a} Barreto, {\em Recovering Asymptotics
of Metrics from Fixed Energy Scattering Data}, to appear in {\it
Invent. Math.}
\bibitem{sslaes} R.B. Melrose, {\em  Spectral and Scattering Theory for
the Laplacian on Asymptotically Euclidean spaces (M. Ikawa, ed),}{  Marcel
Dekker, 1994.}
\bibitem{smagai} R.B. Melrose, M. Zworski, {\em Scattering Metrics and
Geodesic Flow at Infinity,}{ Invent. Math. 124, 389-436 (1996).}
\bibitem{andras} A. Vasy, {\em Geometric Scattering Theory for
Long-Range Potentials and Metrics,} I.M.R.N. 1998 no 6, 285-315
\end{thebibliography}
\end{document}